\documentclass[preprints,commentary,accept,moreauthors,fleqn,pdftex]{Definitions/mdpi} 

\firstpage{1} 
\pubvolume{xx}
\issuenum{1}
\articlenumber{5}
\pubyear{2019}
\copyrightyear{2019}
\history{ }

\pdfoutput=1

\usepackage{latexsym}

\usepackage{amssymb}
\usepackage{amsfonts}

\numberwithin{equation}{section}
\numberwithin{theorem}{section}
\numberwithin{definition}{section}
\numberwithin{property}{section}

\newcommand{\dis}{\displaystyle}

\newcommand{\ba}{\begin{array}}
\newcommand{\ea}{\end{array}}

\newcommand{\be}{\begin{equation}}
\newcommand{\ee}{\end{equation}}

\newcommand{\ods}{\par \vspace{0.2cm} \par}

\newenvironment{my-itemize}{
  \begin{list}{$\bullet$}{
  \setlength{\itemsep}{10pt}
  \setlength{\parskip}{5pt}
  \setlength{\parsep}{15pt}
  \setlength{\leftmargin}{40pt}
  \setlength{\rightmargin}{0pt}}
}{
  \end{list}
}

\newcommand{\bi}{\begin{my-itemize}}
\newcommand{\ei}{\end{my-itemize}}

\Title{{On the famous unpublished } preprint \\
 ``Methods of integration which preserve the contact transformation property
 of the Hamilton equations'' \\ 
by Ren\'e De Vogelaere}

\Author{Robert D. Skeel $^{1}$, \ \    Jan L. Cie\'sli\'nski $^{2}$}

\AuthorNames{Robert D. Skeel and Jan L.\ Cie\'sli\'nski}

\address{%
$^{1}$ \quad School of Mathematical and Statistical Sciences,
Arizona State University, Tempe, AZ 85287, USA; rskeel@asu.edu \\ 
$^{2}$ \quad Uniwersytet w Bia\l ymstoku, Wydzia{\l} Fizyki, ul.\ Cio\l kowskiego 1L,  15-245 Bia\l ystok, Poland;  j.cieslinski@uwb.edu.pl}

\abstract{The well known 1956 unpublished report by De Vogelaere, motivated by the physics of particle accelerators, 
is the first demonstration of the existence of numerical integrators (now called symplectic integrators)  
that preserve a basic property of Hamiltonian systems. The purpose of this short note is to introduce Ren\'e De Vogelaere,
list his pioneering accomplishments,
and relate them to other early work on symplectic integrators. The preprint itself, typset in LaTeX, is attached as an appendix. }

\keyword{De Vogelaere; symplectic integrator; Hamiltonian systems}

\MSC{01A60; 65P10}

\begin{document}



Ren\'e De Vogelaere was born in Etterbeck, Belgium, on May 18, 1926, and died  in Berkeley, USA,  on December 14,  1991. He was a student of Georges Lema\^{i}tre  at the Catholic University of Louvain in Belgium,  receiving a Ph.D. in mathematics in 1948. 
He taught first at Laval University in Quebec,
then at Notre Dame University in Indiana,
where he wrote this well known and farsighted manuscript,
and finally, { since 1958,} at the University of California, Berkeley~\cite{CCLP92}.

The problem of numerically integrating Hamiltonian systems
on very long time intervals was brought to De Vogelaere's attention
by computational physicist J. Snyder.%
\footnote{Snyder
later headed the Computer Science Department at the University of Illinois.}
It was in connection with the design of a particle accelerator
to be built in the Midwestern United States.\footnote{
The project was later cancelled, soon after L.B. Johnson became president.}
De Vogelaere thought that if the numerical integrator
possessed the same qualitative properties as the Hamiltonian system,
then the effect of discretization error would be no worse than
the effect of errors in the model,
a kind of thinking associated with backward error analysis.
The resulting 1956 Notre Dame report~\cite{DeVo56}
produces a number of positive results on numerical symplectic integration.
It shows that the two variants of the method
that later became known as the symplectic Euler method
satisfy the property of being symplectic.
The report characterizes this method as being ``known'',
for the case of a {\em separable} Hamiltonian system.
The report also generalizes the method to a nonseparable system.
By composing the two variants of the nonseparable symplectic Euler method,
De Vogelaere obtains second order schemes for nonseparable systems.
These generalize each of the two variants
of the St\"ormer-Verlet/leapfrog method.
The extension of staggered schemes to nonseparable systems
is unique among those researchers who independently
discover the existence of symplectic numerical integrators.

There is no evidence of any further work on symplectic integration
until the 1983 article by R.D. Ruth~\cite{Ruth83},
who again is motivated by the dynamics of particle accelerators
but unaware of the unpublished report of De Vogelaere.
Using generating functions, Ruth shows that the leapfrog method is symplectic
and derives a novel third order symplectic method,
which employs a Jacobian matrix.
These results are for separable Hamiltonians.

Later that year and again unaware of prior work~\cite{ChSc90},
physicist P.J. Channell writes an internal report~\cite{Chan83}
showing the possibility of doing symplectic integration,
again using generating functions.

A year later, mathematician K.~Feng~\cite{Feng85,Feng86}
independently discovers symplectic numerical integrators,
in particular, the implicit Gauss-Legendre Runge-Kutta methods
and the (staggered) St\"ormer-Verlet method for separable Hamiltonians,
and illustrates the use of generating functions.
Publication was in a mathematics journal that was just three years old,
and the work initially received little attention.

However, by 1988,
a number of mathematical researchers had became
aware of and interested in the development
of symplectic (or canonical) integrators.
In particular, the seminal 1988 paper of Ge \& Marsden~\cite{GeMa88}
 acknowledges De Vogelaere, Channell, and Feng as developers
of the first symplectic integrators. { It is significant that in this paper  Marsden acknowledges also personal interactions with Ren\'e De Vogelaere, in fact they both were at Berkeley at that time.} 
Awareness of symplecticness
by the numerical ordinary differential equation community
became more widespread after the invited paper by J.M. Sanz-Serna
at the IMA Conference on Computational Differential Equations in 1989.
The published version~\cite{Sanz92} cites also the work of Ruth.
For a comprehensive account of subsequent developments,
consult Ref.~\cite{HaLW02}.

The report ``Methods of integration which preserve the contact transformation property  of the Hamilton equations'' \cite{DeVo56}, written in 1956 and still being an unpublished preprint, can be considered as a pioneering work in the field presently known as  geometric numerical integration, far ahead of its time.  {\it  ``This is a marvellous paper, short, clear, elegant, written in one week, submitted for publication -- and never published'' } (E. Hairer, C. Lubich and G. Wanner, see \cite{HaLW03}, p. 406). 

The content of the preprint can be easily summarized in few words, using modern terminology. Ren\'e De Vogelaere has shown that the methods presently known as symplectic Euler methods and St\"ormer-Verlet (or leap-frog) method are symplectic. 

Three cases are considered, in order of increasing complications:
\bi
   \item separable: \  $\ddot x = f (x)$
\item  non-separable: \  $\dot x = f (x,y)$  , \quad $\dot y = g (x,y)$  , \quad  $\displaystyle \frac{\partial f}{\partial y} + \frac{\partial g}{\partial x} = 0$  , 
\item multi-dimensional: \ $\displaystyle   \frac{d q_i}{d t} = \frac{\partial H}{\partial p_i}$  , \quad 
   $\displaystyle  \frac{d p_i}{d t} = - \frac{\partial H}{\partial q_i}$  , \quad 
   $(i = 1,\ldots,n)$  .
 \ei
In the first case notation $\ddot x = g (x)$ would be more convenient, because then this case be obtained from the second one as a special case $f (x,y)=y$, $g(x,y) = g (x)$.   

The symplectic Euler methods (Eqs. (3.1) and (3.3) in \cite{DeVo56}) are the only ones referered to as ``known methods'' by De Vogelaere. He extends the method (3.1) on the second and third case (see Eqs. (5.1) and (7.3) in \cite{DeVo56}) showing that they are symplectic of order 1. The extension of Eq. (3.3) on the second and third case is left as an obvious exercise. Hairer, Lubich and Wanner present these results as Theorem 3.3,  attributing it explicitly to Ren\'e De Vogelaere, see \cite{HaLW02}, p. 189. The role of De Vogelaere in the development of the symplectic Euler methods is discussed also by Blanes and Casas \cite{BC16}, p. 22--23.

Another family of symplectic integrators, leap-frog methods, has even longer history \cite{HaLW03}. In the separable (multidimensional) case they began to be widely used in molecular dynamics only after Verlet's publication in 1967, however they were rediscovered many times before by astronomers and physicists (e.g., by Delambre 1792, Encke 1860 and St\"ormer 1907) and also by\ldots De Vogelaere \cite{DeVo56}. Loup Verlet himself found that this method was in fact used by Isaac Newton (in {\em Principia}, 1687) to prove Kepler's second law \cite{HaLW03}. The surprising effectiveness of this method is due to its symplecticity and this propoerty has been first noticed and proved by De Vogelaere in \cite{DeVo56}. His Eqs. (4.1), (4.3),(6.1) and (7.4) represent, in fact, the leap-frog method:  
\[  \begin{array}{l}
\displaystyle 
\pmb{q}_{1/2} = \pmb{q}_0 + \frac{1}{2} h  \ \frac{\partial H}{\partial \pmb{p}} 
(\pmb{q}_{1/2}, \pmb{p}_0) \ , \qquad 
\pmb{p}_{1/2} = \pmb{p}_0 - \frac{1}{2} h \ \frac{\partial H}{\partial \pmb{q}} (\pmb{q}_{1/2}, \pmb{p}_0)\ , \\[4ex]
\displaystyle
\pmb{q}_1 = \pmb{q}_{1/2} + \frac{1}{2} h  \ \frac{\partial H}{\partial \pmb{p}} 
(\pmb{q}_{1/2}, \pmb{p}_1) \ , \qquad 
\pmb{p}_1 = \pmb{p}_{1/2} - \frac{1}{2} h \ \frac{\partial H}{\partial \pmb{q}} (\pmb{q}_{1/2}, \pmb{p}_1)\ .
\end{array} \]
As already suggested this is the first known construction of the leap-frog method in the case of general non-separable canonical Hamiltonian systems. 

The title of the preprint is a little bit misleading because De Vogelaere considers only transformations in the phase space $(q,p)$ treating time as a parameter. Now, such transformations are rather called symplectic (or canonical). Contact transformations,  usually related to the extended phase space (including the time variable), are more general.

The problem of extending symplectic integrators on non-autonomous (time-dependent) Hamiltonian systems is much more difficult and still only particular 
results are available. It seems that in this aspect the way of introducing time dependence of integrators presented in the famous preprint is not so fruitful as the symplecticity of the integrators in the autonomous case. 

Ren\'e De Vogelaere wrote his report in just one week. James Snyder brought this problem to his attention at the 
MURA\footnote{Midwestern Universities Research Association, a collaboration between 15 universities with the goal of designing and building a particle accelerator.} meeting of April 14, 1956 and the preprint appeared on April 21, 1956.   
Taking into account the surprisingly short time of writing this report, the number of mistakes seems to be very few. Actually, we noticed only one place with some computational errors (not very essential and without any influence on the conclusions of the paper): coefficents by several terms of the operator $\Phi$, appearing in Eq.~(6.2) in section 6, are not correct. It should read as follows:
\[
\begin{array}{l}  \displaystyle
\Phi = - \frac{1}{24} f^2 \frac{\partial^2}{\partial x^2} + \frac{1}{12} g^2 \frac{\partial^2}{\partial y^2} + \left( \frac{\alpha^2 - \alpha}{2} + \frac{1}{12} \right) \frac{\partial^2}{\partial t^2}     + \left(  \frac{1}{6} - \frac{\alpha}{2}  \right) g  \frac{\partial^2}{\partial y \, \partial t} 
\\[3ex] \dis \qquad  - \frac{1}{12} f \frac{\partial^2}{\partial x \, \partial t}  -  \frac{1}{12} f g \frac{\partial^2}{\partial x \, \partial y} 
 + \frac{1}{12} \frac{\partial g}{\partial t} \frac{\partial }{\partial y} + \left( \frac{\alpha}{2} - \frac{1}{6} \right) \frac{\partial f}{\partial t} \frac{\partial }{\partial x} 
 \\[4ex] \displaystyle \qquad
 + \frac{1}{12} f  \frac{\partial g}{\partial x} \frac{\partial }{\partial y} 
 + \frac{1}{12} f  \frac{\partial f}{\partial x} \frac{\partial }{\partial x} + \frac{1}{12} g  \frac{\partial g}{\partial y} \frac{\partial }{\partial y} - \frac{1}{6} g \frac{\partial f}{\partial y}  \frac{\partial }{\partial x} \ . 
\end{array} 
\]

The De Vogelaere's preprint \cite{DeVo56}, typeset in LaTeX for Readers' convenience (the original report is written on a  mechanical typewriter and is hardly readable in some places), is  available as a supplementary material  (appendix)  to our short article.  We tried to preserve the original form of the manuscript, except correcting few obvious misprints and taking also into account the  handwritten corrections made on the preprint, probably by the author  himself. 

\ods\ods



\acknowledgments{We thank Charles DeVogelaere, PhD, and the family of Ren\'e DeVogelaere for granting permission to publish the preprint \cite{DeVo56}. Many thanks are due to the University of Notre Dame and, especially, to J.~Parker Ladwig,  for assistance and help.}

\ods\ods



\reftitle{References}

\end{document}




\section{Introduction}

In recent years the construction of high
powered accelerators has led to vast programs of computations
involving the solution $S_1$  over a very long time range of
Hamiltonian systems describing approximately the motion $S$ 
of a proton in the accelerator. Errors between the computed
solution $S_2$ and $S_1$ are introduced because of the method
of integration (finite step of increase of $t$) and because
of round-off (finite number of binary or decimal digits).
Their worst effect is probably to destroy the contact transformation
property of $S$. Hence the suggestion of using a
method of integration which, if there was no round-off error,
would give a solution $S_3$ with the contact transformation
property. Moreover, if the error, due to the finiteness of
increase of $t_1$ is not too large, one may even expect that
the error between $S_3$ and $S$ will be of the same order of
magnitude as the error between $S_1$ and $S$. How good the
solution $S_2$ has to be, will depend on a study of the
physical system and its Hamilton approximation. One should also make sure that the discontinuities in the derivatives introduced at
each step by the method of integration do not  significantly alter the results.

One may expect that second order methods may lead
to significant results even with step increments which are
not too small.

The contact transformation property does reduce to
area conservation when the Hamiltonian has one degree of
freedom. We feel that it is not superfluous to treat
first this case in detail, it did lead us to the general
case and we hope that the procedure used may give a lead
to the construction of higher order methods. Two known
first order methods are given for the special case
\be    \label{1.1}
  \ddot x = f (x, t) ,
\ee
they are extended to the case
\be \ba{l}  \label{1.2}
\dis \dot x = \frac{\partial H}{\partial y} \ , \qquad 
\dis \dot y = - \frac{\partial H}{\partial x} \ .
\ea \ee 
Second order methods are described. The error term for
one step is given, it may also provide a lead to higher
order methods. The error after $n$ steps is not given.
(7.4) and (7.5) are probably the most important formulae
in this paper.

\section{Some basic properties}

We will first recall a certain
number of basic properties relevant to this paper, although
their generalization is well known.

\begin{Definition}  \label{Def 2.1}
If $x_1 = f (x_0, y_0)$ and $y_1 = g (x_0, y_0)$ are two functions, continuous in  $x_0$, $y_0$ as well as their first partial derivatives,
\[
  \dis J = \frac{ d (x_1, y_1)}{d (x_0, y_0)} = \left|   \ba{cc}  \dis \frac{\partial x_1}{\partial x_0} & \dis \frac{\partial x_1}{\partial y_0} \\[2ex] \dis \frac{\partial y_1}{\partial x_0} & \dis \frac{\partial y_1}{\partial y_0} \ea      \right|
\]
is called the Jacobian of $x_1$, $y_1$ with respect to $x_0$, $y_0$. 
\end{Definition}

\begin{Property}  \label{prop 2.1}
If $x_1 = x_0$ and $y_1 = f (x_0, y_0)$, 
\[
\frac{ d (x_1, y_1)}{d (x_0, y_0)} = \frac{\partial y_1}{\partial y_0} . 
\]
\end{Property}

\begin{Property}  \label{prop 2.2}
If $x_2 = F (x_1, y_1)$ and $y_2 = G (x_1, y_1)$, then
\[
\frac{ d (x_2, y_2)}{d (x_0, y_0)} = \frac{ d (x_2, y_2)}{d (x_1, y_1)} \cdot \frac{ d (x_1, y_1)}{d (x_0, y_0)} . 
\]
\end{Property}

\begin{Property} \label{prop 2.3}
If $\dis \frac{ d (x_1, y_1)}{d (x_0, y_0)} \neq 0$, then 
\[
\frac{ d (x_0, y_0)}{d (x_1, y_1)} = \left( \frac{ d (x_1, y_1)}{d (x_0, y_0)} \right)^{-1} . 
\]
\end{Property}

\begin{Property}  \label{prop 2.4}
The elements of area \ $d x_0  d y_0$ \ and \ $d x_1  d y_1$ \ are related by 
\[
d x_1  d y_1 =  d x_0  d y_0
\]
provided $J=1$.
\end{Property}

\noindent For proofs one may consult de la Vall\'ee Poussin (\cite{vallee}, vol. I, p. 360--364).
\ods

If one considers the differential equations
\be \ba{l}   \label{dotxy}
    \dot x = f (x, y, t) \ , \\[1ex]
    \dot y = g (x, y, t) \ ,  
\ea \ee
where $f$ and $g$ are continuous functions in $x$, $y$ and $t$ and have first partial derivatives continuous in $x$ and $y$; if $x(t)$, $y(t)$ are solutions of \rf{dotxy} with initial conditions $x_0$ and $y_0$ when $t=t_0$; then one has

\begin{Property}  \label{prop 2.5}
\[
   J = \frac{ d (x, y)}{d (x_0, y_0)} = e^{\int_{t_0}^t dt \left( \frac{\partial f}{\partial x} + \frac{\partial g}{\partial y} \right) }
\]
\end{Property} 

\noindent (see for instance de la Vall\'ee Poussin (\cite{vallee}, vol. II, p. 146).

\begin{Property}   \label{prop 2.6}
A necessary and sufficient condition that
\[
   f d y - g d x
\]
be an exact differential $d H$ --when $t$ is considered as a parameter-- is
\[
    \frac{\partial f}{\partial x} + \frac{\partial g}{\partial y}  = 0 ,
\]
 \end{Property} 
 \noindent (\cite{vallee}, vol. II, p. 41).
\ods
Hence

\begin{Theorem}   \label{Th 2.1}
A necessary and sufficient condition that
the equations \rf{1.2} --under the above condition-- be such
that the area in the space of initial conditions at $t = 0$ 
be preserved for all $t$, is that there exists a function
$H(x,y,t)$ such that
\[
  f = \frac{\partial H}{\partial y} \ , \qquad g = - \frac{\partial H}{\partial x} 
\]
or that \rf{dotxy} is a Hamiltonian system. 
\end{Theorem}
\noindent (See also Whittaker \cite{Whitt}, Ch. XI).
\ods

This theorem shows that the only systems of two equations
of first order that we have to consider are of the
form \rf{1.2}. The preceeding properties will enable us to
prove that the given methods of integration have Jacobian
unity, hence are area preserving. These methods of integration
are dependent upon a parameter $h$ (time increment). If
the solution of \rf{dotxy} has a continuous $n$th derivative for
$0 \leqslant  t - t_0 \leqslant h$, one can write
\[
x (t_0 + h) = x(t_0) + h \dot x (t_0) + \ldots + \frac{h^n}{n!} [ x^{(n)}(t_0) + \ep ],
\]
where $\ep$ tends to zero with $h$.

\begin{Definition}  \label{Def 2.2}
If the method of integration is such that
the approximation $\bar x (t_0 + h)$ expanded in terms of $h$ coincides
with the above Taylor expressions to the order $p \leqslant n$, one  
says that the method is of the {\bf $\pmb{p}$th  order} and we will use
 the symbol $\{ p \}$ beside the formulae of the method.
\end{Definition}

We have not included in this paper the algebraic manipulation
which proves rigorously, by this identification
method, that the formulae are correct to the given order. 
We refer to Ince \cite{Ince} or Hildebrandt \cite{Hil} 
for the famous application of this method to obtain
Runge-Kutta type formulae.

If $p \leqslant n-1$, it is possible to give the error term of
the formula as the coefficient of $h^{p+1}$ of $\bar x (t_0 + h)-x(t_0+h)$ multiplied by $h^{p+1}$. This error will be given. The error
after $n$ steps where $nh$ is fixed is of the order of $h^p$.  It
can be found, if necessary, by methods used in similar cases. See for instance Collatz \cite{Col} or Hildebrandt \cite{Hil}.

\section{The known methods}
\label{sec 3}

The following methods are well known
when $\alpha = 0$, but they are only first order and applicable to
the special case \rf{1.1}:
\be  \ba{l}  \label{3.1}
 \bar x_1 = x_0 + h \dot x_0 \\[1ex]
 \dot{\bar x}_1 = \dot x_0 + h f (\bar x_1, t_0 + \alpha h) 
 \ea  \qquad \{ 1 \}
\ee
The error terms are 
\be \ba{l}     \label{3.2}
\dis \bar x_1 - x(t_0 + h) = - \frac{1}{2} h^2 f (x_0, t_0) \\[2ex]
\dis  \dot{\bar x}_1 - \dot x(t_0 + h) = \left[ \frac{1}{2} \dot x_0 \left( \frac{\partial f}{\partial x} \right)_0 + \left(\alpha - \frac{1}{2} \right) \left( \frac{\partial f}{\partial t} \right)_0 \,    \right] h^2 ,
\ea \ee
the last equation suggests to take $\alpha = \frac{1}{2}$. Also
\be \ba{l}    \label{3.3}
\dot{\bar x}_1 = \dot x_0 + h f (x_0, t_0 + \alpha h) \\[1ex]
\bar{x}_1 = x_0 + h \dot x_1 
\ea \qquad \{ 1 \}
\ee
the error terms are
\be \ba{l}    \label{3.4}
\dis \bar x_1 - x(t_0 + h) =  \frac{1}{2} h^2 f (x_0, t_0) \\[2ex]
\dis  \dot{\bar x}_1 - \dot x(t_0 + h) = \left[ - \frac{1}{2} \dot x_0 \left( \frac{\partial f}{\partial x} \right)_0 + \left(\alpha - \frac{1}{2}\right) \left( \frac{\partial f}{\partial t} \right)_0  \,   \right] h^2 ,
\ea \ee
the last equation suggests to take $\alpha = \frac{1}{2}$.

It is easy to check the area conservation property
using definition~\ref{Def 2.1}.

Alternately, one may also write the transformation
as a product of two transformations and use the property~\ref{prop 2.2}.
For instance~\rf{3.1}  can be written
\[  \ba{l}
\bar x_1 = x_0 + h \dot x_0 \ , \qquad  \dot x_0 = \dot x_0 \\[1ex]
\dot{\bar x}_1 = \dot x_0 + h f ({\bar x}_1, t_0 + \alpha h) \ , \qquad {\bar x}_1 = {\bar x}_1 \ ;
\ea \]
the Jacobian of both transformations is one because of the property~\ref{prop 2.1}.

\section{Second order methods for~\rf{1.1}}
\label{sec 4}

The following method can be used to solve~\rf{1.1}:
\be  \ba{l}  \label{4.1}
\dot {\bar x}_1 = \dot x_0 + \frac{1}{2} h f (x_0, t_0 + \alpha h) \\[1ex]
\bar{x}_2 = x_0 + h \dot {\bar x}_1  \\[1ex]
\dot {\bar x}_2 = \dot {\bar x}_1  + \frac{1}{2} h f (\bar{x}_2, t_0 + (1-\alpha) h) 
\ea \qquad \{ 2 \}
\ee

One can infer that the method is second order because
the second equation uses the slope at the mid-point, and the
last equation combined with the first uses the mean of the
slopes of $\dot x$ at the extremities of the interval, indeed the 
error terms are
\be \ba{l}  \label{4.2}
{\bar x}_2 - x (t_0 + h) = \dis \left[ - \frac{1}{6} \dot x \, \frac{\partial f}{\partial x} + \left( \frac{\alpha}{2} - \frac{1}{6} \right)  \frac{\partial f}{\partial t}    \right]_0 h^3 \\[4ex] \dis
\dot {\bar x}_2 - \dot x (t_0 + h) = \left[  \frac{1}{12} \frac{\partial^2 f}{\partial x^2} \, {\dot x}^2 + \left( \frac{1}{6} - \frac{\alpha}{2}  \right)
\frac{\partial^2 f}{\partial x \partial t} \, \dot x +  \left( \frac{\alpha^2 - \alpha}{2} + \frac{1}{12}  \right) \frac{\partial^2 f}{\partial t^2} + \frac{1}{12} f \frac{\partial f}{\partial x} 
\right]_0 h^3
\ea 
\ee

The alternate method at the end of section~\ref{sec 3}, gives
immediately the area conservation property. The first
equation suggests that we take $\alpha = \frac{1}{3}$. 

A similar argument shows that we can also use
\be  \ba{l}    \label{4.3}
{\bar x}_1 =  x_0 + \frac{1}{2} h  {\dot x}_0  \\[1ex]
\dot {\bar x}_2 = \dot x_0 + h f (\bar x_1, t_0 + \frac{1}{2} h)  \\[1ex]
{\bar x}_2 =  \bar{x}_1  + \frac{1}{2} h \dot {\bar x}_2 
\ea \qquad \{ 2 \}
\ee
with the error terms:
\be \ba{l}  \label{4.4}
{\bar x}_2 - x (t_0 + h) = \dis \left[ \frac{\partial f}{\partial x} \, \dot x +  \frac{\partial f}{\partial t}  \right]_0 \frac{h^3}{12}  \\[4ex] \dis
\dot {\bar x}_2 - \dot x (t_0 + h) = - \left[  \frac{1}{2} \frac{\partial^2 f}{\partial x^2} \, {\dot x}^2 + 
\frac{\partial^2 f}{\partial x \partial t} \, \dot x + 
 \frac{1}{2} \frac{\partial^2 f}{\partial t^2} + 2 \frac{\partial f}{\partial x} \, f 
\right]_0 \frac{h^3}{12} 
\ea 
\ee

This second method has comparable error terms, but asks
for only one computntion of $f( x, t)$ per interval, hence is
barely more complicated that the first order method.

\section{First order method for \rf{1.2} }
\label{sec 5}

The methods \rf{3.1} and \rf{3.3} are symmetrical to each
other and do suggest the following generalization:
\be  \ba{l}  \label{5.1}
{\bar x}_1 =  x_0 + h f (\bar x_1, y_0,  t_0 + \alpha h)   \\[1ex]
{\bar y}_1 =  y_0 + h g (\bar x_1, y_0,  t_0 + \alpha h)
\ea \qquad \{ 1 \}
\ee
with the error term
\[  \ba{l} \dis
{\bar x}_1 - x(t_0+h) = \left[ - \frac{1}{2} \frac{\partial}{\partial y} (f g) + (\alpha - \frac{1}{2} )  \frac{\partial f}{\partial t}   \right]_0 h^2  \\[3ex]  \dis
{\bar y}_1 - y(t_0+h) = \left[  \frac{1}{2} \frac{\partial}{\partial x}  (f g) + (\alpha - \frac{1}{2} ) \frac{\partial g}{\partial t}   \right]_0 h^2 
\ea \] 
We suggest $\alpha = \frac{1}{2}$.

The easiest way to prove that the Jacobian is one is
to write the transformation as a product of transformations
and to use the properties~\ref{prop 2.2}, \ref{prop 2.3} (for the  first transformation) and~\ref{prop 2.6}.

The first equation~\rf{5.1} gives $x_1$  by an implicit formula,
hence the method is lengthier than the special case
\rf{3.1} obtained when $f(x,y,t) \equiv  y$  or than the other special  case when $f(x,y,t)$ is linear in $x$.  In all other cases, the solution will be obtained by iteration; because $y$ and $t$ are fixed any accelerative process of iteration will furnish
quickly the solution. Tho most obvious accelerating
processes are Newton's and Aitken's methods.

The simple iteration method defined by the equation
\rf{5.1} will converge if

\be
\left|  h \frac{\partial f}{\partial x}  (x, y_0, t_0 + \alpha h) \right| < 1
\ee
Because $h$ must be taken small (to make the error terms
small enough), this relation will not usually lead to an
additional restriction and the rate of convergence will
usually be fast.

It is clear that one may interchange $x$ and $y$ in \rf{5.1},
the form \rf{3.3} is then a special case.

\section{Second order method for \rf{1.2}}
\label{sec 6}

The methods \rf{4.1} and \rf{4.3} did suggest to us the
generalization
\be  \ba{l}   \label{6.1} 
{\bar x}_1 = x_0 + \frac{h}{2} f ({\bar x}_1, y_0, t_0 + \alpha h) \\[2ex]
{\bar y}_1 = y_0 + \frac{h}{2} g ({\bar x}_1, y_0, t_0 + \alpha h) \\[2ex]
{\bar y}_2 = {\bar y}_1 + \frac{h}{2} g ({\bar x}_1, {\bar y}_2, t_0 + (1-\alpha) h) \\[2ex]
{\bar x}_2 = {\bar x}_1 + \frac{h}{2} f ({\bar x}_1, {\bar y}_2, t_0 + (1-\alpha) h)
\ea \qquad \{ 2 \}
\ee 
with the error terms
\be \ba{l}     \label{6.2}
{\bar x}_2 - x(t_0 + h) = (\Phi f)_0 \ h^3 \\[2ex]
{\bar y}_2 - y(t_0 + h) = (\Phi g)_0 \ h^3 
\ea \ee
where the operator $\Phi$ is given by
\[  \ba{l}  \dis
\Phi = - \frac{1}{24} f^2 \frac{\partial^2}{\partial x^2} + \frac{1}{12} g^2 \frac{\partial^2}{\partial y^2} + \left( \frac{\alpha^2 - \alpha}{2} + \frac{1}{12} \right) \frac{\partial^2}{\partial t^2}     + \left(  \frac{1}{6} - \frac{\alpha}{2}  \right) g  \frac{\partial^2}{\partial y \, \partial t} - \frac{1}{12} f \frac{\partial^2}{\partial x \, \partial t}  -  \frac{1}{12} f g \frac{\partial^2}{\partial x \, \partial y} 
\\[3ex] \dis
\qquad + \frac{1}{4} \frac{\partial g}{\partial t} \frac{\partial }{\partial y} + \frac{\alpha}{2}  \frac{\partial f}{\partial t} \frac{\partial }{\partial x} + \frac{1}{4} f  \frac{\partial g}{\partial x} \frac{\partial }{\partial y} 
 + \frac{1}{2} f  \frac{\partial f}{\partial x} \frac{\partial }{\partial x} + \frac{1}{4} g  \frac{\partial g}{\partial y} \frac{\partial }{\partial y} - \frac{1}{6} g \frac{\partial f}{\partial y}  \frac{\partial }{\partial x} \ . 
\ea \]
We suggest $\alpha=0$.

The two first equations are analogous to \rf{5.1} hence
the Jacobian of ${\bar x}_1, {\bar y}_1$ with respect to $x_0, y_0$ is one. The
two last equations are analogous to \rf{5.1}  where $x$ and $y$ 
are interchanged, hence tho Jacobian of ${\bar x}_2, {\bar y}_2$ with respect
to ${\bar x}_1, {\bar y}_1$ is one. The complete transformation preserves the 
area because of the property~\ref{prop 2.2}.

That the method is second order could have been inferred
by remarking that
\[
{\bar y}_2 = y_0 + \frac{h}{2} \ g ({\bar x}_1, y_0, t_0 + \alpha h) + 
\frac{h}{2} \ g ({\bar x}_1, y_2, t_0 + (1-\alpha) h) 
\]
hence that if $y$ did not appear in $g$ the slope is taken at
the mid-point for $x$ and that if $x$ did not appearing one 
uses the mean of the slopes at the extremity of the interval
for $y$.

Two of tho relations \rf{6.1}  are implicit, this seems to
be the price we have to pay for the general equations, when
one insists on a method which preserves the area.

Of course, one may interchange $x$ and $y$ in \rf{6.1}.

\section{Generalization to $n$ degrees of freedom}

Let us now consider the Hamilton system
\be   \ba{l}
\dis q'_i = \frac{\partial H}{\partial p_i}  \\[3ex] \dis
p'_i = - \frac{\partial H}{\partial q_i} 
\ea   \qquad H=H(q_i, p_i, t) \ ,  \quad  i=1,2,\ldots, n 
\ee  
and the solutions $q_i (\alpha_i, \beta_i, t)$, $p_i (\alpha_i, \beta_i, t)$ with initial conditions  $\alpha_i, \beta_i$ when $t = t_0$. 

The property $J=1$ generalizes into Poisson brackets relations (see~\cite{Whitt}, Ch. XI):
\be \ba{l}  \label{7.2}
[p_i, p_k] = 0 \\[2ex]
[q_i, q_k] = 0 \\[2ex]
[p_i, q_k] = \delta_{ik} \ .
\ea   \qquad i,k = 1,2,\ldots,n 
\ee 
To check that a transformation is a contact one, one may
either use \rf{7.2} or prove that
\[
  \sum (p_i d q_i - \beta_i d \alpha_i )
\]
is a total differential, when $p$ and $q$ are expressed in
terms of $\alpha$  and $\beta$ and when $t$ is considered as a parameter.
The method
\be   \ba{l}  \label{7.3}
\dis q_i = \alpha_i + h \frac{\partial H (q_i, \beta_i, t)}{\partial \beta_i}  \\[4ex] \dis
p_i = \beta_i - h \frac{\partial H (q_i, \beta_i, t)}{\partial q_i} 
\ea   \qquad \{ 1 \} 
\ee  
is first order. See \rf{5.1} which suggests $t=t_0 + \frac{1}{2} h$. \rf{7.3} is a contact transformation:
\[  \ba{l} \dis
\sum (p_i d q_i - \beta_i d \alpha_i ) = \sum \left( p_i d q_i + \alpha_i d \beta_i - d (\alpha_i \beta_i) \right)  \\[3ex] \dis
\quad \qquad = \sum \left( \left( \beta_i - h \frac{\partial H}{\partial q_i} \right) d q_i + \left( q_i - h \frac{\partial H}{\partial \beta_i} \right) d \beta_i - d (\alpha_i \beta_i)  \right) \\[4ex] \dis
\quad \qquad = d \left(   \sum  \left( \beta_i q_i  - \beta_i \alpha_i - H (q_i, \beta_i, t) \right)    \right) 
\ea \]
where $t$ is considered as a parameter.

Similarly the method 
\be \ba{l}  \label{7.4}
\dis q_i = \alpha_i + \frac{h}{2}  \frac{\partial H (q_i, \beta_i, t_0 + \alpha h)}{\partial \beta_i}  \\[4ex] \dis
p_i = \beta_i - \frac{h}{2} \frac{\partial H (q_i, \beta_i, t_0 + \alpha h)}{\partial q_i}    \\[4ex] \dis
p_i (t_0 + h) \cong P_i = p_i - \frac{h}{2} \frac{\partial H (q_i, P_i, t_0 + (1-\alpha) h)}{\partial q_i}  \\[4ex] \dis
q_i (t_0 + h) \cong Q_i = q_i + \frac{h}{2} \frac{\partial H (q_i, P_i, t_0 + (1-\alpha) h)}{\partial P_i}
\ea \qquad \{ 2 \} 
\ee
is second order as seen in section~\ref{sec 6}. \rf{7.4} is a contact
transformation because it is the product of two contact
transformations. We suggest $\alpha=0$.

In general two of the above relations are implicit.
The special case
\[
   H = \frac{1}{2} \left( \sum p_i^2  + 2 U (q_i, t)  \right)
\]
is worth mentioning. \rf{7.4} reduces then to
\be \ba{l}  \label{7.5}
\dis q_{i1} = q_{i 0}  + \frac{1}{2} h\ {\dot q}_{i 0}  \\[3ex] \dis
 q_i (t_0 + h) \cong q_{i 2} = {\dot q}_{i 0}  - h \ \frac{\partial U (q_{i 1}, t_0 + \frac{1}{2} h)}{\partial q_{i 1}}  \\[3ex] \dis
q_i (t_0 + h) \cong q_{i 2} = q_{i 1}  + \frac{1}{2} h\ {\dot q}_{i 2} .
\ea \qquad \{ 2 \} 
\ee

No detailed example or discussion is given. This will
best be done by those working on these problems in the
Brookhaven, Harwell, MURA or CERN group.
\ods\ods

\noindent {\bf 8. Acknowledgment}.  I wish to thank Professor J. Snyder of the University of Illinois and MURA for having brought this problem to my attention at the MURA meeting of April 14, 1956.

\ods\ods





\reftitle{References}

\ods\ods \ods\ods
\ods\ods \ods\ods
\ods\ods \ods\ods
\ods\ods \ods\ods
\ods\ods \ods\ods
\noindent {\bf Editor's note:} 

The above text is a faithful representation (in terms of content, not necessarily the layout) of the original preprint, i.e., the Report 4 of Center of Numerical Analysis, University of Notre Dame, Indiana, \url{https://curate.nd.edu/downloads/6682x34921d}, except its title page.  
The original form is preserved except correcting few obvious misprints and taking into account handwritten notes made on the manuscript, probably by the author himself. Some informations from the title page (like affiliation) are included, but we ommited ``Internal Classification: me 2.4, pr 5.2, pr 5.4'', ``under contract with the Office of Naval Research N7-0nr-43906 (NR-041-148)'', and ``submitted for publication''.

We thank Charles DeVogelaere, PhD, and the family of Ren\'e DeVogelaere for granting permission to publish this paper. Many thanks are due to the University of Notre Dame and, especially, to J.~Parker Ladwig,  for assistance and help. 

\hfill J. L. Cie\'sli\'nski